\documentclass[amssymb,11pt,reqno]{amsart}

\usepackage{amsfonts}
\usepackage{amssymb}
\usepackage{amsthm}
\usepackage{newlfont}

\newcommand{\R}{\mathbb{R}}
\newcommand{\C}{\mathbb{C}}

\newcommand{\N}{\mathbb{N}}
\newcommand{\G}{\Gamma}
\newcommand{\h}{\mathcal H}

\newcommand{\test}{\mathcal S}

\def\pf{\noindent {\bf Proof :  \  }}
\def\vol{\mbox{\rm Vol}}
\def\a{\alpha}

\def\sg{\sigma}
\def\rkn{\R^{\k n}}
\def\skn{S^{\k n-1}}
\def\rsg{R_{\sg}}

\def\k{\kappa}
\def\hxip{\h_{\xi}^{\bot}}
\def\hxi{\h_{\xi}}
\newtheorem{thm}{Theorem}

\newtheorem{lemma}{Lemma}
\newtheorem{prop}{Proposition}

\usepackage{graphicx}
\usepackage{latexsym}
\usepackage{enumerate}
\begin{document}

\title{A note on the Busemann-Petty problem for bodies of certain invariance}

\author{Marisa Zymonopoulou}
\address{
Marisa Zymonopoulou\\
Department of Mathematics\\
Case Western Reserve University \\
Cleveland, OH 44106, USA}

\email{marisa@@cwru.edu}

\begin{abstract}The Busemann-Petty problem asks whether origin symmetric convex bodies in $\R^n$ with smaller hyperplane sections necessarily have smaller volume. The answer is affirmative if $n\leq 3$ and negative if $n\geq 4.$ We consider a class of convex bodies that have a certain invariance property with respect to their ordered $\k$-tuples of coordinates in $\rkn$ and prove the corresponding problem.
\end{abstract}
  
\maketitle

\section{Introduction}\label{intro}

In 1956 the Busemann-Petty problem was posed (see [BP]), asking the following question:

\noindent
Suppose that $K$ and $L$ are two origin symmetric convex bodies in $\R^n$ such that for every $\xi \in S^{n-1},$

$$\vol_{n-1}\bigl(K\cap \xi^{\perp}\bigr) \leq \vol_{n-1}\bigl(L\cap \xi^{\perp}\bigr).$$
Does it follow that
$$\vol_{n}\bigl(K\bigr) \leq \vol_{n}\bigl(L \bigr) \ \ ?$$

\noindent
The answer is affirmative if $n\leq 4$ and negative if $n\geq 5.$ The problem was solved in the late 90's as a result of a series of papers ( [Ba], [Bo], [Ga1], [Ga2], [GKS], [Gi], [H], [K1], [K2], [LR], [Lu], [Pa],  [Zh1],  [Zh2], ; see [K8, p.3] for the history of the solution).

In this article we consider the corresponding problem in $\rkn$ for convex bodies that have the property of invariance with respect to certain rotations considering the volumes of the $\k n-\k$-dimensional central sections, where $\k\in\N.$

Let $x=(x_1,x_2,\ldots , x_{\k n})\in\rkn.$ For every $\sg \in SO(\k)$ we define 
$$\rsg(x):=\bigl(\sg(x_1,\ldots,x_{\k}),\cdots, \sg(x_{\k(n-1)+1},\ldots,x_{\k n})\bigr)$$
to be the vector that corresponds to the rotation of the ordered $\k $-tuples of coordinates of $x$ by $\sg.$  

We consider the class of convex bodies that are invariant under any $R_\sg$ rotations:
 \noindent
 if $x=(x_1,x_2,\ldots , x_{kn})\in\rkn$, and $\sg\in SO(\k),$ then 
$$\|x\|_D=\|\bigl(\sg(x_1,\ldots,x_{\k}),\cdots, \sg(x_{\k(n-1)+1},\ldots,x_{\k n})\bigr) \|_D.$$
We then say that the body $D$ has the \emph{$\rsg$-invariance property} or that it is \emph{$\rsg$-invariant}.

Now suppose $\{I_d, \sg _1, \ldots, \sg _{\k-1}\}\subset SO(\k)$ is an orthogonal sequence of counterclockwise rotations. Then the set $\{I_d, R_{\sg _1},\ldots,R_{\sg _{\k-1}} \}$ is also orthogonal, namely, for every $x \in\rkn$ the vectors  $ x, R_{\sg _1}(x),\ldots,R_{\sg _{\k-1}}(x)$ are mutually orthogonal.
 For $\xi \in \skn,$ we denote by $\hxip$ the $\k$-dimensional subspace of $\rkn$ generated by $\{R_j(\xi)\}_{j=0}^{k-1},$ (where $R_0=I_d).$

\medskip

Now, the problem can be formulated as follows: 

\noindent
Suppose $K$ and $L$ are two origin symmetric, with the $\rsg$-invariance property, convex bodies in $\rkn,$ so that
$$\vol_{\k n-\k} (K\cap H_\xi) \leq \vol_{\k n-\k} (L\cap H_\xi),$$
for every $\xi \in \skn.$ Does it follow that
$$\vol_{\k n}(K)\leq \vol_{\k n}(L) \ ?$$

Note that when $\k=1$ the problem corresponds to the real Busemann-Petty problem and for $\k=2$ the problem coincides with the complex Busemann-Petty problem, which has affirmative answer if $n\leq 3$ and negative if $n\geq 4.$ 

In this article we prove that the answer to the Busemann-Petty problem for bodies with the $\rsg$-invariance property is affirmative in the cases where $(i)\ n=2, \ \k\in\N, \ (ii)\ n=3 ,\ \k\leq 2 ,\ (iii)\ n=4$ and $\k=1$ and negative for any other values of $n$ and $\k.$

 A related result on the Busemann-Petty problem for bodies that satisfy an invariance property can be found in [Ru].
\bigbreak

\section{The analytic connection to the problem}\label{appr}

Throughout this paper we use the Fourier transform of distributions. The Schwartz class of rapidly decreasing infinitely differentiable functions (test functions) in $\R^n$ is denoted by $\test(\R^n),$ and the space of distributions over $\test(\R^n)$ by $\test^{\prime}(\R^n).$ The Fourier transform $
at{f}$ of a distribution $f \in \test^{\prime}(\R^n)$ is defined by $\langle\hat{f},\phi\rangle=\langle f,\hat{\phi}\rangle$ for every test function $\phi.$ A distribution is called even homogeneous of degree $p \in \R$ if $\langle f(x),\phi(x/\alpha)\rangle=|\alpha|^{n+p}\langle f,\phi\rangle$ for every test function $\phi$ and every $\alpha \in \R, \ \alpha\neq 0.$ The Fourier transform of an even homogeneous distribution of degree $p$ is an even homogeneous distribution of degree $-n-p.$ A distribution $f$ is called positive definite if, for every test function $\phi, \ \langle f, \phi\ast\overline{\phi(-x)}\rangle \geq 0.$ By Schwartz's generalization of Bochner's theorem, this is equivalent to $\hat{f}$ being a positive distribution in the sense that  $\langle\hat{f},\phi\rangle \geq 0$ for every non-negative test function $\phi,$ (see [K8, section 2.5] for more details).

We denote by $\Delta$ the Laplace operator and by $|\cdot|_2$ the Euclidean norm in the proper space.

 A compact set $K \subset\R^n$ is called a star body, if every straight line that passes through the origin crosses the boundary of the set at exactly two points and the boundary of $K$ is continuous in the sense that the \emph{Minkowski functional} of $K,$ defined by
$$\|x\|_K=\min \{\alpha \geq 0 : x \in \alpha K \},$$
is a continuous function on $\R^n.$
Using polar coordinates it is possible to obtain the following \emph{polar formula for the volume} of the body:

$$\vol_n(K)=\int_{\R^n}\chi (\|x\|_K)dx=\frac{1}{n}\int_{S^{n-1}}\|\theta\|_K^{-n}d\theta.$$

A star body $K$ in $\R^n$ is called $k$-smooth (infinitely smooth) if the restriction of $\|x\|_{K}$ to the sphere $S^{n-1}$ belongs to the class of $C^k(S^{n-1} )\  (C^{\infty}(S^{n-1})).$ It is well-known that one can approximate any convex body in $\R^n$ in the radial metric,
$d(K, L)=\sup \{|\rho_{K}(\xi)-\rho_{L}(\xi)|,\ \xi \in S^{n-1} \},$
by a sequence of infinitely smooth convex bodies. The proof is based on a simple convolution argument (see for example [Sch, Theorem 3.3.1]). The same idea can be used to show that any convex body in $\rkn$ with the $\rsg$-invariance property can be approximated in the radial metric by a sequence of infinitely smooth convex bodies that are also $\rsg$-invariant. This follows from the fact that the $\rsg$-invariance is preserved under convolutions.

If $D$ is an infinitely smooth origin symmetric star body in $\R^n$ and $0<k <n,$ then the Fourier transform of the distribution $\|x\|_D^{-k}$ is a homogeneous function of degree $-n+k$ on $\R^n,$ whose restriction to the sphere is infinitely smooth (see [K8, Lemma 3.16]).

The following Proposition is a spherical version of Parseval's formula established in [K4], (see also [K8, Lemma 3.22]):
\begin{prop}\label{parseval}
 Let $K$ and $L$ be two infinitely smooth origin symmetric convex bodies in $\R^n$ and $0<p<n.$ Then

 $$\int_{S^{n-1}}\bigl(\|x\|_K^{-p}\bigr)^{\wedge}(\xi)\bigl(\|x\|_L^{-n+p}\bigr)^{\wedge}(\xi)d\xi
 =(2\pi)^n\int_{S^{n-1}}\|x\|_K^{-p}\|x\|_L^{-n+p}dx.$$
\end{prop}
\medskip

The class of intersection bodies was introduced by Lutwak [Lu]. This class was generalized in [K3], as follows:
Let $1\leq k <n,$ and let $D$ and $L$ be two origin symmetric star bodies in $\R^n.$ We say that $D$ is the \emph{$k$-intersection body of $L$} if for every $(n-k)-$dimensional subspace $H$ of $\R^n$
$$\vol_k(D\cap H^{\perp})=\vol_{n-k}(L\cap H ).$$
More generally, we say that an origin symmetric star body $D$ in $\R^n$ is a {\it $k$-intersection body}
if there exists a finite Borel measure $\mu$ on $S^{n-1}$ so that for every even test function $\phi\in {\cal{S}}(\R^n),$
$$\int_{\R^n} \|x\|_D^{-k} \phi(x)\ dx =
\int_{S^{n-1}} \left(
\int_0^\infty t^{k-1} \hat\phi(t\xi)\ dt\right) d\mu(\xi).$$
Note that $k$-intersection bodies of star bodies are those $k$-intersection
bodies for which the measure $\mu$ has a continuous strictly positive density; see [K7] or [K8, p. 77].

A more general concept of embedding in $L_{-p}$ was introduced in [K6].
Let $D$ be an origin symmetric star body in $\R^n,$ and
$X=(\R^n,\|\cdot\|_D).$ For $0<p<n,$
we say that $X$ embeds in $L_{-p}$
if there exists a finite Borel measure $\mu $ on $S^{n-1}$ so that, for every even
test function $\phi$
$$
\int_{\R^n} \|x\|_D^{-p} \phi(x)\ dx = \int_{S^{n-1}} \left(
\int_{\R} |z|^{p-1} \hat{\phi}(z\theta)\ dz\right) d\mu(\theta).
$$
Obviously, an origin symmetric star body $D$ in $\R^n$ is a $k$-intersection body
if and only if the space $(\R^n,\|\cdot\|_D)$ embeds in $L_{-k}.$ In this article we use
embeddings in $L_{-p}$ only to state some results in continuous form;  for more
applications of this concept, see [K8, Ch. 6].

Embeddings in $L_{-p}$ and $k$-intersection bodies admit a Fourier analytic
characterization  that we are going to use throughout
this article:\begin{prop} \label{posdef} ([K7], [K8, Th. 6.16]) Let $D$ be an origin symmetric star body in $\R^n,\ 0<p<n.$ The space
$(\R^n,\|\cdot\|_D)$ embeds in $L_{-p}$ if and only if the function $\|x\|_D^{-p}$ represents
a positive definite distribution on $\R^n.$ In particular, $D$ is a $k$-intersection body if
and only if $\|x\|_D^{-k}$ is a positive definite distribution on $\R^n.$
\end{prop}
\medskip

Another important fact that will be used is the co-called second derivative test, (see [K8, Theorems 4.19, 4.21]).
\begin{prop} \label{sdt} Let $\lambda\ge 3,\ k\in {\mathbb N}\cup \{0\}, \ q>2$ and
let $Y$ be a finite dimensional normed
space of dimension greater or equal to $\lambda.$ Then the $q$-sum of $\R$ and $Y$
does not embed in $L_{-p}$ with $0<p<\lambda-2.$ In particular, this direct sum is not
a $k$-intersection body for any $1\le k < \lambda-2.$
\end{prop}
Recall that the $q$-sum of $\R$ with a space $Y,$ $Y\oplus_q \R,$ 
is defined as the space of pairs $\{(t,y):\ \in \R, y\in Y\}$ with the norm
$\|(y,t)\|= \left(\|y\|_Y^q+t^q \right)^{1/q}.$

\medskip

\bigskip

Let $ \k\geq 1$ and let $H$ be an $(\k n-\k)$-dimensional subspace of $\R^{\k n}.$
Fix any orthonormal basis $e_1,...,e_\k$  in the orthogonal subspace $H^\bot.$
For a convex body $D$ in $\rkn,$ define the $(\k n-\k)$-dimensional parallel section
function $A_{D,H}$ as a  function on $\R^\k$ such that
$$A_{D,H}(u) = \vol_{\k n-\k}(D\cap \{H+ u_1e_1+...+u_\k e_\k\})$$
 \begin{equation} \label{paral}
 = \int\limits_{\{x\in \rkn: (x,e_1)=u_1,...,
(x,e_k)=u_\k\}} \chi(\|x\|_D)\ dx, \quad u\in \R^\k.
\end{equation}

If $D$ is infinitely smooth the function $A_{D,H}$ is infinitely differentiable at the origin (see [K8, Lemma 2.4]). So we can consider the action of the distribution $|u|_2^{-q-\k}/ \G(-q/2)$ on $A_{D,H}$ and apply a standard regularization argument (see for example [K8, p.36] and [GS, p.10]). Then the function
\begin{equation} \label{fraclapl}
q\mapsto \Big\langle {{|u|_2^{-q-\k}}\over {\Gamma(-q/2)}},
A_{D,H}(u) \Big\rangle
\end{equation}
is an entire function of $q\in \C.$ 
If  $q=2m,\ m\in \N\cup \{0\},$ then
$$ \Big\langle {{|u|_2^{-q-\k}}\over {\Gamma(-q/2)}}\Big\vert_{q=2m},
A_{D,H}(u) \Big\rangle$$
\begin{equation} \label{intaction}
 = \frac{(-1)^m |S^{\k-1}|}{2^{m+1} \k(\k+2)...(\k+2m-2)} \Delta^m A_{D,H}(0),
\end{equation}
where $|S^{\k-1}| = 2\pi^{\k/2}/\Gamma(\k/2)$ is the surface area of
the unit sphere $S^{\k-1}$ in $\R^\k.$ When the body $D$ is
origin symmetric, the function $A_{D,H}$ is even, and for $0<q<2$
we have (see also [K8 p. 39])
$$\Big\langle \frac{|u|_2^{-q-\k}}{\Gamma(-q/2)}, A_{D,H}(u) \Big\rangle $$
\begin{equation} \label{action}
= \frac{1}{\Gamma(-q/2)}\int_{S^{\k n-1}} \left( \int_0^\infty \frac{A_{D,H}(t\theta) -
A_{D,H}(0)}{t^{1+q}}\ dt \right)\ d\theta.
\end{equation}
Note that the function (\ref{fraclapl}) is equal (up to a constant) to the fractional power
of the Laplacian $\Delta^{q/2}A_{D,H},$(see [KKZ] or [K7] for complete definition).

\medskip

The following proposition was proved in [KKZ, Prop. 4].

\begin{prop} \label{mainprop} Let $D$ be an infinitely smooth origin symmetric
convex body in $\rkn$ and $1\le \k <\k n.$ Then for every $(\k n-\k)$-dimensional
subspace $H$ of $\rkn$ and any $q\in \R,\ -\k<q< \k n-\k,$
$$\Big\langle {{|u|_2^{-q-\k}}\over {\Gamma(-q/2)}},
A_{D,H}(u) \Big\rangle$$
\begin{equation} \label{eq1}
 = {{2^{-q-\k}\pi^{-\k/2}}\over {\Gamma((q+\k)/2)(\k n-q-\k)}}
\int_{S^{\k n-1}\cap H^\bot} \big(\|x\|_D^{-\k n+q+\k}\big)^\wedge(\theta)
\ d\theta.
\end{equation}
Also for every $m\in {\mathbb{N}} \cup \{0\},\ m<(\k n-\k)/2,$
\begin{equation} \label{eq2}
 \Delta^m A_{D,H}(0) = {{(-1)^m}\over {(2\pi)^k(\k n-2m-\k)}}
\int_{S^{\k n-1}\cap H^{\bot}} \bigl(\|x\|_D^{-\k n+2m+\k}\bigr)^\wedge(\eta)\ d\eta.
\end{equation}

\end{prop}

\medbreak

Brunn's theorem (see for example [K8, Theorem 2.3]) states that for an origin symmetric convex body and a fixed direction, the central hyperplane section has the maximal volume among all the hyperplane sections perpendicular to the given direction. As a consequence we have the following generalization proved in [KKZ, Lemma 1].
\begin{lemma} \label{brunn} If $D$ is a $2$-smooth origin symmetric convex body in $\rkn,$
then the function $A_{D,H}$ is twice differentiable at the origin and
$$\Delta A_{D,H}(0) \le 0.$$ Besides that for any $q\in (0,2),$
$$\Big\langle \frac{|u|_2^{-q-\k}}{\Gamma(-q/2)}, A_{D,H}(u) \Big\rangle \ge 0.$$
\end{lemma}

\bigskip

                          \section{Connection with $\k $-intersection bodies}

Now, we are ready to connect the $\rsg$-invariance property of the bodies with the Fourier trasform of their norm that will yield the solution to the problem. The following simple observation is crucial for applying the analytic methods to convex bodies that carry the $\rsg$-invariance property.

\begin{lemma} \label{const} Suppose that $D$ is an origin symmetric infinitely smooth,
with the $R_\sg$-invariance property, star body in $\rkn.$ Then for
every $0<\a<\k n$ and $\xi\in \skn$ the Fourier transform of the distribution $\|x\|_D^{-\a}$ is
a constant function on $\skn\cap \hxip.$
\end{lemma}
\pf By [K8, Lemma 3.16], since $D$ is infinitely smooth,   $\bigl(\|x\|_D^{-\a}\bigr)^{\wedge}$ is a continuous
function outside of the origin in $\rkn.$  In addition, since the body $D$ has the $\rsg$-invariance property, by the
connection between the Fourier transform of distributions and linear transformations,
the Fourier transform of $\|x\|_D^{-\a}$ is also $\rsg$-invariant. 
Now, as mentioned in the Introduction, the $\k$-dimensional subspace $\hxip$ of $\rkn$ is spanned by 
$\xi, R_{\sg_1}(\xi), \ldots, R_{\sg_{\k-1}}(\xi)$. Hence, every vector in the $(\k-1)$-dimensional sphere $\skn\cap \hxip$
is the image of $\xi$ under one of the coordinate-wise rotations $R_\sg$, so
the Fourier transform of $\|x\|_D^{-\a}$ is a constant function on $\skn\cap \hxip.$
\qed

\medskip

We use the above lemma to express the volume
of $\k$-codimensional sections in terms of the Fourier transform.

\begin{thm} \label{volconst} Let $K$ be an infinitely smooth origin symmetric, with the $\rsg$-invariance property, convex
body in $\rkn, n\ge 2.$  For every $\xi\in \skn,$ we have
$$\vol_{\k n- \k}(K\cap \hxi) = \frac{2^{-\k+1}\pi^{-\frac{\k}{2}}}{(\k n-\k)\G(\k/2)} \left(\|x\|_K^{-\k n+\k}\right)^\wedge(\xi).$$
\end{thm}

\pf
Let $\xi\in\skn.$ In formula (\ref{eq2}) we put $H=\hxi$ and $m=0.$ Then, by the definition of the $(\k n-\k)$-dimensional section function (equation (\ref{paral})), we have that
$$A_{D,\hxi}(0)=\vol_{\k n-\k}(K\cap \hxi)=\frac{1}{(2\pi)^{\k}(\k n-\k)}\int_{\skn\cap \hxip}
\left(\|x\|_K^{-\k n+\k}\right)^\wedge(\eta)d\eta.$$
By Lemma \ref{const}, the function under the integral is constant on the $(\k-1)$-dimensional sphere $\skn\cap \hxip.$ Since $\xi\in\skn\cap \hxip,$  we get that
$$ \vol_{\k n-\k}(K\cap \hxi)=\frac{|S^{\k-1}|}{(2\pi)^{\k}(\k n-\k)\G(\k/2)}\left(\|x\|_K^{-\k n+\k}\right)^\wedge(\xi),$$
which proves the theorem, using the volume of the unit sphere of $\R^{\k}.$
\qed
\bigskip

As in the case of the real Busemann-Petty (for $\k=1$) and the complex Busemann-Petty (where $\k=2$), the solution to this problem is closely related to the $\k$-intersection bodies. 

\begin{thm} \label{connection} The answer to the Busemann-Petty problem for convex bodies with the $\rsg$-invariance property
is affirmative if and only if every origin symmetric $\rsg$-invariant 
convex body in $\rkn$ is a $\k$-intersection body.
\end{thm}

The proof of the theorem follows by the next two lemmas. Note that we may assume that the bodies are infinitely smooth. This follows from the fact that one can approximate the two bodies, in the radial metric, with infinitely smooth convex bodies that have the $\rsg$-invariance property (see Section \ref{appr}). Then, if the answer to the problem is affirmative for infinitely smooth bodies it is also true for the original bodies.

\medskip

\begin{lemma}\label{lm:pos}
Suppose $K$ and $L$ are infinitely smooth origin symmetric, with the $R_\sg$-invariance property, convex bodies in $\R^{2n}$ so that $K$ is a $\k$-intersection body. Then, if for every $\xi \in \skn$
\begin{equation}\label{sectineq}
\vol_{\k n-\k} (K\cap \hxi) \leq \vol_{\k n-\k} (L\cap \hxi),\end{equation}
then
 $$\vol_{\k n}(K)\leq \vol_{\k n}(L). $$

\end{lemma}

\pf
The bodies $K$ and $L$ are infinitely smooth, so by [K8, Th. 3.16], the Fourier transform of the distributions $\|\cdot\|_K^{-\k n+\k}, \ \|\cdot\|_K^{-\k}$ and $\|\cdot\|_L^{-\k n+\k}$ are continuous functions outside of the origin of $\rkn.$ By Theorem \ref{volconst}, the inequality (\ref{sectineq}) becomes
$$ \left(\|x\|_K^{-\k n+\k}\right)^\wedge(\xi) \le \left(\|x\|_L^{-\k n+\k}\right)^\wedge(\xi).$$
Also, since $K$ is a $\k$-intersection body, Proposition
\ref{posdef} implies that the Fourier transform of the distribution $\|\cdot\|_K^{-\k}$ is a non-negative function. So, we can multiply both sides of the latter inequality by $\bigl(\|\cdot\|_K^{-\k}\bigr)^{\wedge}$ . We then integrate over $\skn$ and use Parseval's spherical version, Proposition \ref{parseval}, to have that
$$ \int_{\skn}\bigl(\|\cdot\|_K^{-\k}\bigr)^{\wedge}(\xi)\left(\|x\|_K^{-\k n+\k}\right)^\wedge(\xi)d\xi$$
$$ \le \int_{\skn}\bigl(\|\cdot\|_K^{-\k}\bigr)^{\wedge}(\xi)\left(\|x\|_L^{-\k n+\k}\right)^\wedge(\xi)d\xi $$
which gives
$$ \int_{\skn}\|x\|_K^{-\k n}dx \le \int_{\skn}\|x\|_K^{-\k}\|x\|_L^{-\k n+\k}dx. $$
We use the polar formula for the volume and H\"older's inequality, then the latter inequality becomes
$$\k n \vol_{\k n}(K)\leq \Bigl(\int_{\skn}\|x\|_K^{-\k n}dx\Bigr)^{\frac{1}{n}}\Bigl(\int_{\skn}\|x\|_L^{-\k n}dx\Bigr)^{\frac{n-1}{n}},$$
which proves the Lemma.
\qed

\medskip

For the negative part we need a perturbation argument to construct a body that will give a counter-example to the problem. The new body immediately inherits the additional property of invariance with respect to the $\rsg$ rotations of the original convex body.

\begin{lemma}\label{lm:neg}
Suppose that there exists an origin symmetric  convex body $L$ in $\rkn,$ with the $R_\sg$-invariance property,
which is not a $\k$-intersection body. Then there exists an origin symmetric $R_\sg$-invariant convex body $K$ in $\rkn$ such that for every $\xi\in \skn,$
$$\vol_{\k n-\k}(K\cap \hxi) \le \vol_{\k n-\k}(L \cap \hxi),$$
but
$$\vol_{\k n}(K) > \vol_{\k n}(L).$$
\end{lemma}

\pf
We assume that the body $L$ is infinitely smooth with positive curvature. By [K8, Lemma 3.16] the Fourier transform of the distribution $\|x\|_{L}^{-\k}$ is a continuous function on the unit sphere $\skn.$ Moreover there exists an open subset $\Omega$ of $\skn$ in which $\bigl(\|x\|_{L}^{-\k}\bigr)^{\wedge}<0 .$
Since $L$ is invariant with respect to all $R_\sg$ rotations we may assume that $\Omega$ has also this invariance property.
 
We use a standard perturbation procedure for convex bodies, see for example [K8, p.96] (similar argument was used in [KKZ, Lemma 5]):
 Consider a non-negative infinitely differentiable even function $g$ supported on $\Omega$ that is also $R_\sg$-invariant.  We extend it to a homogeneous function of degree $-\k$ on $\R^{\k n}.$ By [K8, Lemma 3.16] its Fourier transform is an even homogeneous function of degree $-\k n+\k$ on $\R^{\k n},$ whose restriction to the sphere is infinitely smooth: $\bigl(g(x/|x|_2)|x|_2^{-\k}\bigr)^{\wedge}(y)=h(y/|y|_2)|y|_2^{-\k n+\k},$ where $h\in C^{\infty}(\skn).$

We define a body $K$ so that
$$\|x\|_{K}^{-\k n+\k}=\|x\|_{L}^{-\k n+\k}-\varepsilon |x|_2^{-\k n+\k}h\bigl(\frac{x}{|x|_2}\bigr),$$
for small enough $\varepsilon >0$ so that the body $K$ is strictly convex. Note that $K$ also has the $\rsg$-invariance property. We apply the Fourier transform to both sides and the latter becomes
\begin{eqnarray}\label{fcountex}
\bigl(\|x\|_{K}^{-\k n+\k} \bigr)^{\wedge}(\xi)&=&\bigl(\|x\|_{L}^{-\k n+\k} \bigr)^{\wedge}(\xi)
- (2\pi)^{\k n}\epsilon g(\xi) \\
&\leq & \left(\|x\|_L^{-\k n+\k}\right)^{\wedge}(\xi),\nonumber\end{eqnarray}
since $g$ is non-negative. By Theorem \ref{volconst}, this implies that for every $\xi\in \skn$
$$\vol_{\k n-\k}(K\cap H_\xi) \le \vol_{\k n-\k}(L \cap H_\xi).$$
Now, we multilpy both sides of (\ref{fcountex}) by $\bigl(\|\cdot\|^{-\k}\bigr)^{\wedge}$ and integrate over the $\skn.$ Then
$$\int_{\skn} \left(\|x\|_K^{-\k n+\k}\right)^\wedge(\xi) \left(\|x\|_L^{-\k}\right)^\wedge(\xi)\  d\xi$$
$$ = \int_{\skn} \left(\|x\|_L^{-\k n+\k}\right)^\wedge(\xi) \left(\|x\|_L^{-\k}\right)^\wedge(\xi)\  d\xi$$
$$- (2\pi)^{\k n}\epsilon \int_{\skn} \left(\|x\|_L^{-\k}\right)^\wedge(\xi) g(\xi)\ d\xi$$
$$> \int_{\skn} \left(\|x\|_L^{-\k n+\k}\right)^\wedge(\xi) \left(\|x\|_L^{-\k}\right)^\wedge(\xi)\  d\xi,$$
since $\left(\|x\|_L^{-\k}\right)^\wedge<0$ on the support of $g.$
 We apply Proposition \ref{parseval} and the polar formula for the volume to obtain
$$\int_{\skn} \|x\|_K^{-\k n+\k}\|x\|_L^{-\k}  d\xi > \k n \vol_{\k n}(L).$$

As in Lemma \ref{lm:pos}, we use H\"older's inequality and the polar formula for the volume of $K$ to obtain the inverse inequality for the $\k n$-dimensional volumes of the bodies.
\qed

\section{The solution}

In order to prove the main result we need to determine the dimensions and the values of $\k$ for which every origin symmetric, with the $\rsg$-invariance property, convex body in $\rkn$ is a $\k$-intersection body.
It is known (see [K5] or [K8, Corollary 4.9] with Proposition \ref{posdef}) that for every origin symmetric convex body $D$ in $\rkn, \  n\geq3,$ the space $(\rkn,\|\cdot\|_D)$ embeds in $L_{-p}$ for $p\in [kn-3,kn).$ In other words, every origin symmetric convex body $D$ in $\rkn,$ is a $(\k n-3),\ 
(\k n-2), \ (\k n-1)$-intersection body.

\medskip

\noindent
\textbf{Remark: }Let $n\geq 3.$ We consider the unit ball
$$
B_q^{\k n} = \{x\in \rkn:
 \|x\|_q=\bigl((x_1^2+\cdots +x_\k^2)^{\frac{q}{2}}+...$$
$$\cdots+ (x_{\k (n-1)+1}^2+\cdots +x_{\k n}^2)^{\frac{q}{2}}\bigr)^{1/q}\le 1\},
$$
with $q>2.$ The space $\bigl(\rkn, \|\cdot\|_{B_q^{\k n}}\bigr)$ contains the $q$-sum of $\R$ with the $(\k n-\k)$-dimensional subspace, $Y=\bigl(\R^{\k n-\k},\|\cdot\|_{B_q^{\k n-\k}}\bigr).$ Then by
the second derivative test, Proposition \ref{sdt}, for $\lambda=\k n-\k,$ the space $\R\oplus_q Y$ does not embed in $L_{-p},$ if $p<\k n-\k-2.$ Moreover, it follows by a result of E.Milman [Mi] that the larger space cannot embed in $L_{-p},$ for $p<\k n-\k-2.$ Hence, if $\k>\frac{2}{n-2},$ the unit ball $B_q^{\k n}$ provides a counter-example of an origin symmetric, with the $\rsg$-invariance property, convex body that is not a $\k$-intersection body.

\medskip

It remains to see what happens in the range $[\k n-\k -2, \k n-3).$

\begin{thm}\label{kn-k-2} Every origin symmetric, with the $\rsg$-invariance property, convex body $K$ in $\rkn$ is a $(\k n-\k-2)$-intersection body in the following cases: $(i)\ n=2, \k\geq3,\ (ii)\ n=3, \k\geq 2, (iii)\ n\geq 4, \k\in\N.$

Moreover, the space $(\rkn, \|\cdot\|_K)$ embeds in $L_{-p}$ for $p\in [\k n-\k-2,\k n),$ when $(i)\  n=2, \k\geq2, \ (ii)\ n\geq 3, \k\in\N.$


\end{thm}

\pf
We use Proposition \ref{mainprop} with $m=1$ and $H=\hxi.$ The condition $1<\frac{\k n-\k}{2}.$ is valid for the three cases of the first part of the theorem.
Hence, by (\ref{eq2}) and Lemma \ref{const}, we have that

$$\Delta A_{D,\hxi}(0) = {(-1)|S^{\k -1|}\over {(2\pi)^k(\k n-2-\k)}}
 \bigl(\|x\|_D^{-\k n+2+\k}\bigr)^\wedge(\xi).$$
By the generalization of Brunn's theorem (Lemma \ref{brunn}), we conclude that 
$$\bigl(\|x\|_D^{-\k n+\k+2}\bigr)^\wedge(\xi)\geq 0, $$ 
for every $\xi\in\skn.$ Hence, by Proposition \ref{posdef}, $K$ is a
$(\k n-\k-2)$-intersection body for the above cases.

Now let $ n=2, \ \k\geq2, $ or $n\geq 3$ with $ \k\in\N.$ By (\ref{eq1}) with $H=\hxi$ and Lemma \ref{const}, we have that for every $q\in (0,2)$ 
$$\Big\langle {{|u|_2^{-q-\k}}\over {\Gamma(-q/2)}},
A_{K,\hxi}(u) \Big\rangle
= \frac{2^{-q-\k}\pi^{-\frac{\k}{2}}|S^{\k -1}|}{\Gamma((q+\k)/2)(\k n-q-\k)}
 \big(\|x\|_K^{-\k n+q+\k}\big)^\wedge(\xi).$$
 Then Lemma \ref{brunn} implies that for $0<q<2,$ 
$\big(\|x\|_K^{-\k n+q+\k}\big)^\wedge\geq 0.$ The range of $q$ and Proposition \ref{posdef} indicate that the
space $(\rkn, \|\cdot\|_K)$ embeds in $L_{-p}$ for $p\in [\k n-\k-2,\k n).$
\qed

\medskip

Using the above theorem we can now prove the main result of this article.

\begin{thm}\label{res}
The Busemann-Petty problem for origin symmetric convex bodies that have the $\rsg$-invariance property in $\rkn $ has affirmative answer only if: $(i)\ n=2, \ \k\in\N, \ (ii)\ n=3 ,\ \k\leq 2 ,\ (iii)\ n=4$ and $\k=1.$
\end{thm}

\pf

When $n=2$ and $\k\in\N$ the answer is affirmative by Theorems \ref{connection}, \ref{kn-k-2} and Proposition \ref{posdef}, since $\k=2\k -\k>2\k-\k-2.$  

If $n=3$ and $\k\leq 2$ then $\k\geq 3\k -\k-2$ and using Theorems \ref{connection}, \ref{kn-k-2} and Proposition \ref{posdef} we have that the answer is also affirmative. On the other hand, if $\k\geq 3$ then $\k<3\k -\k-2.$ By the Theorem \ref{connection} and the Remark the latter implies that the answer to the problem is negative and the body $B_q^{3\k}$ gives a counter-example of a body that is not $\k$-intersection body in $\R^{3\k}.$ 

For $n=4$ and $\k=1$ then $1=\k n-\k-2.$ So the answer is also affirmative. 

Lastly, if $n=4$ with $k\geq 2$ or $n\geq 5$ and $\k\in\N$ then $\k n-\k-2>\k.$ The same reasoning as before, applying Theorem \ref{connection} and the Remark, gives that the answer to the Busemann-Petty problem for origin symmetric convex bodies that have the $\rsg$-invariance property is negative and  $B^{\k n}_q$ provides a counter-example.

\qed

\end{document}